# On the Marginal Distributions of Stationary AR(1) Sequences

**S Satheesh**, Neelolpalam, Trichur - 680 004, **India**.

*ssatheesh1963@yahoo.co.in*

**E Sandhya**, Prajyoti Niketan College, Pudukad, Trichur - 680 301, **India**.

*esandhya@hotmail.com*

**Abstract.** In this note we correct an omission in our paper (Satheesh and Sandhya, 2005) in defining semi-selfdecomposable laws and also show with examples that the marginal distributions of a stationary AR(1) process need not even be infinitely divisible.

**Keywords and Phrases.** Semi-selfdecomposability, infinite divisibility, AR(1) process.

**1. Introduction and Discussion.** The purpose of this note is two-fold. One is to correct an omission in our paper (Satheesh and Sandhya, 2005) in defining semi-selfdecomposable (SSD) laws and (ii) to show with examples that the marginal distributions of an AR(1) process need not even be infinitely divisible (ID). While quoting the definition of SSD laws in terms of characteristic functions (CF) the authors did not mention that the CF $f_o$ should be ID. We correct this omission and discuss the effect of this omission on other results. The correct definition is:

**Definition.1** (Maejima and Naito, 1998). A probability distribution on **R** with CF $f$ is SSD($b$) *if* there exists a CF $f_o$ that is ID such that

$$f(t) = f(bt) f_o(t), \text{ for some } 0<b<1, \forall\, t\in \mathbf{R}. \tag{1}$$

The linear additive AR(1) scheme considered here is described by the sequence of r.vs $\{X_n\}$ and innovations $\{\varepsilon_n\}$ where $\{\varepsilon_n\}$ are *i.i.d*, satisfying

$$X_n = \rho X_{n-1}+ \varepsilon_n, \text{ for some } 0<\rho<1. \tag{2}$$

The innovations $\{\varepsilon_n\}$ in the AR(1) model need not necessarily be ID. Consequently theorem.1 cannot be an "*iff*" statement and it should be corrected as:

**Theorem.1** A sequence $\{X_n\}$ of *r.v*s defines an AR(1) sequence that is marginally stationary with $0<\rho<1$ *if* $X_n$ is SSD($\rho$).

This flaw was then carried over to the description of discrete SSD laws (theorem.9) and the analogous (false) characterization of integer-valued AR(1) sequence in theorem.11. They should be corrected as:

**Theorem.9** An integer-valued distribution on $\{0,1,2,\ldots\}$ with probability generating function (PGF) $P(s)$ is discrete SSD($c$) *if* there exists a PGF $P_o(s)$ that is ID such that

$$P(s) = P(1-c+cs)\, P_o(s), \text{ for some } 0<c<1,\ \forall\, s\in(0,1). \tag{3}$$

Notice that the proof of theorem.9 is still true. In terms of PGFs, a marginally stationary integer-valued AR(1) sequence $\{X_n\}$ with innovations $\{\varepsilon_n\}$ are represented by;

$$P(s) = P(1-\rho+\rho s)\, P_\varepsilon(s), \text{ for some } 0<\rho<1. \tag{4}$$

Here also the innovations $\{\varepsilon_n\}$ need not be ID. Consequently theorem.11 becomes;

**Theorem.11.** A sequence $\{X_n\}$ of integer-valued r.vs defines a marginally stationary AR(1) sequence with $0<\rho<1$ *if* $X_n$ is discrete SSD($\rho$).

It may also be noted that this omission affects only theorems 1 and 11. Other results continue to be true since now we will be considering mixtures of semi-stable or stable laws (which are ID) and where the mixing distribution is also ID and hence the mixture is also ID, Feller (1971, *p*.573).

We now give two examples that illustrate that the expression (not the description) for SSD laws in terms of PGFs (CFs) can be satisfied by a PGF (CF) that is not even ID.

**Example.1** Christoph and Schrieber (2000) showed that the Sibuya($\alpha,\lambda$) law with PGF $P_\lambda(s) = 1-\lambda(1-s)^\alpha$, $0<\alpha\leq 1$, $0<\lambda\leq 1$, is ID *iff* $\lambda\leq 1-\alpha$ and further

$$Q(s) = \frac{1-(1-s)^\alpha}{1-\lambda(1-s)^\alpha}, \ 0<\lambda\leq\frac{1-\alpha}{1+\alpha}, \text{ is also a PGF.}$$

In other words, for some $\lambda \in (0,1)$ the PGF of the Sibuya($\alpha$,1) law satisfies the expression

$$1-(1-s)^\alpha = (1-\lambda(1-s)^\alpha) Q(s) \text{ or}$$

$$P(s) = P(1-\lambda^{1/\alpha} + \lambda^{1/\alpha} s) Q(s), \ \forall \ s \in (0,1),$$

which is similar to (3) with $c = \lambda^{1/\alpha}$ or (4) with $\rho = \lambda^{1/\alpha}$. Also Sibuya($\alpha$,1) law is not ID.

**Example.2** For $\beta>0$ the CF of the triangular($\beta$) law on **R** (Feller, 1971, p.503) is:

$$\varphi_\beta(t) = 1-\beta|t|, \text{ for } |t|\leq \tfrac{1}{\beta} \text{ and } \varphi_\beta(t) = 0, \text{ for } |t|>\tfrac{1}{\beta}.$$

Taking the Sibuya($\alpha,\lambda$)-sum of triangular(1) law shows that $\varphi(t) = 1-\lambda|t|^\alpha$ is also a CF.

Now consider the distribution with CF $Q(\varphi_1(t))$, $Q(s)$ being the PGF in example.1. Then;

$$Q(\varphi_1(t)) = \begin{cases} \dfrac{1-|t|^\alpha}{1-\lambda|t|^\alpha}, & |t|\leq 1 \text{ and} \\ 0, & |t|\geq 1. \end{cases}$$

In other words, for some $\lambda \in (0,1)$, the CF $\varphi(t) = 1-|t|^\alpha$, $0<\alpha\leq 1$ satisfies the relation

$$1-|t|^\alpha = (1-\lambda|t|^\alpha)Q(\varphi_1(t)) \text{ or}$$

$$\varphi(t) = \varphi(\lambda^{1/\alpha} t) \, \varphi_o(t), \ \forall \ t \in \mathbf{R},$$

which is similar to the expression in (1) with $b = \lambda^{1/\alpha}$ or (2) that is marginally stationary in terms of CFs with $\rho = \lambda^{1/\alpha}$. Since the CF $1-|t|^\alpha$ has real zeroes it is not ID.

**Remark.1** Thus the expression (not the description or definition) for SSD laws in terms of CFs (PGFs) can be satisfied by a CF (PGF) that is not even ID. But SSD laws are ID, Maejima and Naito (1998). Hence the condition that the CF $f_o$ is ID in definition.1 and that the PGF $P_o(s)$ is ID in theorem.9 are necessary.

**Remark.2** These examples also show that the marginal distributions of marginally stationary AR(1) sequence need not be SSD or even ID and this conclusion is true in the integer-valued case as well.

Authors sincerely regret that though inadvertent, they wrongly quoted the definition of SSD laws in their work, resulting in incorrect claims.

**Acknowledgement.** The authors sincerely thank Professor N Bouzar, University of Indianapolis, USA, for bringing to their notice the omission and also the reference Christoph and Schrieber (2000).

**References.**